\documentclass[11pt,reqno]{amsart}
\usepackage[utf8]{inputenc}

\usepackage{amssymb, amsmath, amsthm}
\usepackage{hyperref}
\usepackage[alphabetic,lite]{amsrefs}

\usepackage{graphicx}
\usepackage[cmtip,all]{xy}
\usepackage{verbatim}
\usepackage{youngtab}
\usepackage{ytableau}
\usepackage{ upgreek }
\usepackage{ mathrsfs }
\usepackage{verbatimbox}
\usepackage{tikz}
\usepackage{rank-2-roots}

\newtheorem{theorem}{Theorem}[section]

\theoremstyle{definition}
\newtheorem{definition}[theorem]{Definition}

\newtheorem{problem}[theorem]{Problem}
\newtheorem{question}[theorem]{Question}
\newtheorem{conjecture}[theorem]{Conjecture}

\theoremstyle{remark}

\numberwithin{equation}{section}

\newcommand{\defi}[1]{{\bf\upshape\sffamily #1}}
\DeclareMathOperator{\ShHom}{\mathscr{H}\text{\kern -3pt {\calligra\large om}}\,}

\renewcommand{\a}{\alpha}
\newcommand{\bw}{\bigwedge}
\def\kk{{\mathbf k}}
\renewcommand{\ll}{\lambda}
\newcommand{\onto}{\twoheadrightarrow}
\newcommand{\oo}{\otimes}
\newcommand{\pd}{\partial}

\newcommand{\GL}{\operatorname{GL}}

\newcommand{\SL}{\operatorname{SL}}

\newcommand{\Sym}{\operatorname{Sym}}

\newcommand{\bb}[1]{\mathbb{#1}}

\newcommand{\mc}[1]{\mathcal{#1}}
\newcommand{\mf}[1]{\mathfrak{#1}}
\newcommand{\ol}[1]{\overline{#1}}
\newcommand{\op}[1]{\operatorname{#1}}

\newcommand{\ul}[1]{\underline{#1}}

\def\PP{{\mathbf P}}
\def\lra{\longrightarrow}

\begin{document}

\title[Questions about cohomology on flag varieties]{Some questions arising from the study of \\ cohomology on flag varieties}

\address{Department of Mathematics,
University of Notre Dame \hfill \newline
\indent 255 Hurley Notre Dame, IN 46556, USA\newline}

\author{Zhao Gao}
\email{zgao1@nd.edu}

\author{Claudiu Raicu}
\email{craicu@nd.edu}
\thanks{Raicu was supported by the NSF Grant No.~2302341}

\author{Keller VandeBogert}
\email{kvandebo@nd.edu}
\thanks{VandeBogert was supported by the NSF Grant No.~2202871}


\date{}

\begin{abstract} A fundamental problem at the confluence of algebraic geometry and representation theory is to describe the cohomology of line bundles on flag varieties over a field of characteristic~$p$. When $p=0$, the solution is given by the celebrated Borel--Weil--Bott Theorem, while for $p>0$ the problem is widely open. In this note we describe a collection of open questions that arise from the study of particular cases of the general theory, focusing on their combinatorial and commutative algebra aspects.  
\end{abstract}

\maketitle

\section{Introduction}\label{sec:intro}

Let $\kk$ be an algebraically closed field of characteristic $p$, and consider the \defi{(complete) flag variety} $Fl_n$ parametrizing flags of subspaces
\[
 V_{\bullet}:\qquad 0\subset V_1\subset\cdots\subset V_{n-1} \subset \kk^n,\qquad\text{ where }\dim(V_i)=i.
\]
Flag varieties have a rich theory that overlaps many important areas in mathematics: algebraic geometry, representation theory, algebraic group theory, singularity theory, combinatorics, commutative algebra etc. This brief note is motivated by the following foundational question about flag varieties.

\begin{problem}\label{prob:BWB-char-p}
 Describe the cohomology of line bundles on $Fl_n$.
\end{problem}

In characteristic $p=0$, the answer comes from the Borel--Weil--Bott theorem \cite{Serre,Bott} (see also \cite[Section~II.5]{Jantzen}, \cite[Chapter~4]{Weyman}, \cite[Section~11.1]{BCRV}), but in positive characteristic the problem remains open, despite important progress over the years. Perhaps a more modest goal, which remains equally elusive, is to determine which cohomology groups are zero, and which ones are not.

\begin{problem}\label{prob:char-p-vanishing}
 Characterize the vanishing of cohomology of line bundles on~$Fl_n$.
\end{problem}

The best results concerning Problem~\ref{prob:char-p-vanishing} are Kempf's Vanishing theorem \cite{Kempf1,Kempf2}, and Andersen's characterization of the non-vanishing of $H^1$ \cite{Andersen-H1}. Essentially complete answers to Problems \ref{prob:BWB-char-p} and \ref{prob:char-p-vanishing} are known for $n=2$ (when $Fl_n=\PP^1$) and for $n=3$ \cite{Griffith,Donkin,Liu,Gao-Raicu}, but the problems remain open already when $n=4$.

Problems \ref{prob:BWB-char-p} and \ref{prob:char-p-vanishing} can be formulated within the broad context of \defi{generalized flag varieties} $G/B$, where $G$ is a reductive algebraic group, and $B$ is a Borel subgroup, and there is a wealth of important results in this setting that we will not touch upon (see the influential book by Jantzen \cite{Jantzen} and the recent survey by Andersen \cite{Andersen-survey}). The varieties $Fl_n$ correspond to the case when $G=\GL_n$ and $B$ is the Borel subgroup of upper triangular matrices.

Our own interest in Problems~\ref{prob:BWB-char-p} and \ref{prob:char-p-vanishing} comes from the many applications Borel--Weil--Bott theory has had in describing homological invariants in commutative algebra: syzygies \cite{Lascoux,Weyman}, Ext modules \cite{Rai-Ext}, and local cohomology \cite{RWW,RW,LR} for determinantal varieties, syzygies of Segre varieties \cite{ORS}, vanishing results for Koszul modules \cite{AFPRW1,AFPRW2} etc. The main obstacle to extending such results to positive characteristic comes from the difficulty of computing cohomology in positive characteristic, as illustrated in what follows.

\section{Arithmetic complexes}\label{sec:arithmetic}

\tikzset{
  treenode/.style = {align=center, inner sep=0pt, text centered,solid,thin,
    font=\sffamily},
  arn_n/.style = {treenode, circle, white, font=\sffamily\bfseries, draw=black,
    fill=black, text width=.5em},
  arn_nl/.style = {treenode, circle, white, font=\sffamily\bfseries, draw=black,
    fill=black, text width=1.5em},  
  arn_r/.style = {treenode, circle, red, draw=red, 
    text width=.5em, very thick},
  arn_v/.style = {treenode, circle, black, font=\sffamily\bfseries, draw=black, text width=1.2em},
  arn_x/.style = {treenode, rectangle, draw=black,
    minimum width=.5em, minimum height=0.5em},
  dott/.style={edge from parent/.style={dotted, very thick,circle,draw}},
  emph/.style={edge from parent/.style={dashed, very thick,circle,draw}},
  norm/.style={edge from parent/.style={solid,thin,circle,draw}}
}
\makeatletter
\def\labelbox#1{%
  \hbox{%
    \setbox\z@=\hbox{$\m@th\labelstyle{\,#1\,}$}%
    \setbox\tw@=\hbox{$\m@th\labelstyle\,$}%
    \dimen@=\ht\z@ \advance\dimen@ by \wd\tw@ \ht\z@=\dimen@
    \dimen@=\dp\z@ \advance\dimen@ by \wd\tw@ \dp\z@=\dimen@
    \box\z@
  }%
}
\makeatother

Consider $d\geq 1$ and a collection of non-negative integers $w_0,\cdots,w_d$. Write $[d]=\{1,\cdots,d\}$ and think of $[d]$ as indexing edges, and of $w_i$ as weights on the vertices of the path graph
\[ 
\xymatrix{ \overset{w_0}{\bullet} \ar@{-}[r]^1 & \overset{w_1}{\bullet} \ar@{-}[r]^2 & \overset{w_2}{\bullet} \ar@{.}[r] & \bullet \ar@{-}[r]^d & \overset{w_d}{\bullet} }
\]
Every subset $J\subseteq[d]$ determines a decomposition of the path graph into disjoint intervals, and we define the \defi{weight} of an interval to be the sum of weights of the vertices it contains. The removal of an element $j\in J$ breaks up exactly one interval, say of weight $w=w(J,j)$, into two intervals of weights $w'=w'(J,j)$ and $w''=w''(J,j)$, with $w'+w''=w$. We define the complex of $\kk$-vector spaces $C_{\bullet}=C_{\bullet}(w_0,\cdots,w_d)$, where
\[ C_k = \bigoplus_{|J|=k} \kk\cdot e_J\]
and differentials $\pd:C_k\lra C_{k-1}$ given by
\begin{equation}\label{eq:def-pdeJ}
 \pd(e_J) = \sum_{j\in J} (-1)^{s(J,j)} {w(J,j) \choose w'(J,j)}\cdot e_{J\setminus\{j\}},
\end{equation}
where $s(J,j)$ denotes the cardinality of the set $\{i\in[d]:i<j,\ i\not\in J\}$. For an example of a binomial coefficient arising in \eqref{eq:def-pdeJ}, take $J=\{1,2,3\}$, $j=2$:
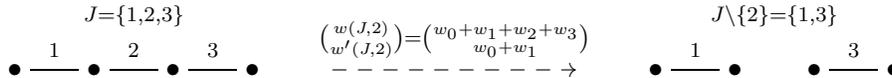
\[
\xymatrix @R=.25pc@C=1.5pc{ 
& \ar@{}[r]^{J = \{1,2,3\}} & & & & & & & & & \ar@{}[r]^{J \setminus\{2\} = \{1,3\}} & & \\
\bullet \ar@{-}[r]^1 & \bullet \ar@{-}[r]^2 & \bullet \ar@{-}[r]^3 & \bullet & \ar@{-->}[rrrr]^{{w(J,2)\choose w'(J,2)}={w_0+w_1+w_2+w_3\choose w_0+w_1}} & & & & & \bullet \ar@{-}[r]^1 & \bullet & \bullet \ar@{-}[r]^3 & \bullet \\
}
\]
The complex $C_{\bullet}$ for $d=3$ is illustrated below in Figure~\ref{fig:Cw-complex}. In the special case when $w_0=w_1=w_2=w_3=1$, we obtain
\begin{equation}\label{eq:C1111}
\xymatrix{0 \ar[r] & \kk \ar[rr]^{\begin{bmatrix} 4 \\ 6 \\ 4 \end{bmatrix}} & & \kk^3 \ar[rrr]^{\begin{bmatrix} -3 & 2 & 0 \\ -3 & 0 & 3 \\ 0 & -2 & 3 \end{bmatrix}}  & & & \kk^3 \ar[rr]^{\begin{bmatrix} 2 & -2 & 2 \end{bmatrix}} & & \kk \ar[r] & 0}
\end{equation}
which is exact if $p=\op{char}(\kk)$ satisfies $p=0$ or $p\geq 5$, but has non-trivial homology when $p=2,3$. If we make the convention that ${0\choose 0}=1$ and take $w_0=\cdots=w_d=0$, then $C_{\bullet}$ is essentially the complex computing the reduced homology of a simplex, which is exact in all characteristics. If all the binomial coefficients ${w\choose w'}$ appearing in the definition of $C_{\bullet}(w_0,\cdots,w_d)$ are non-zero in $\kk$ (for instance if $p=0$ or $p>w_0+\cdots+w_d$), then it is also easy to see that the resulting complex is exact. In general however, the calculation of homology and its dependence on $\kk$ is an open problem.

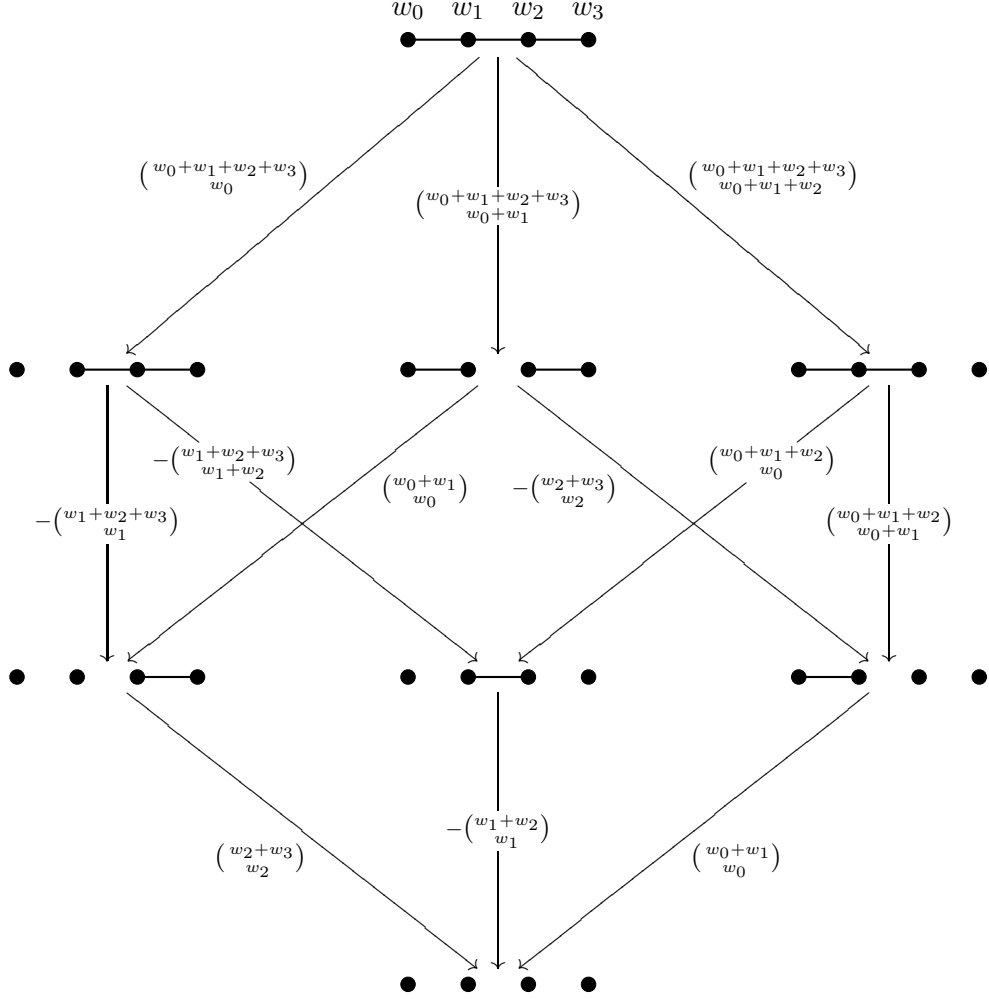
\begin{figure}[h]
\[
\xymatrix @C=5.5em @R=9.5em{
 & {\begin{tikzpicture}[scale=0.8]
\draw [thick] (0,0)--(1,0)--(2,0)--(3,0);
\node[arn_n] at (0,0) {};
\node[label=above:{$w_0$}] at (0,0) {};
\node[arn_n] at (1,0) {}; 
\node[label=above:{$w_1$}] at (1,0) {};
\node[arn_n] at (2,0) {}; 
\node[label=above:{$w_2$}] at (2,0) {};
\node[arn_n] at (3,0) {}; 
\node[label=above:{$w_3$}] at (3,0) {};
\end{tikzpicture}} 
\ar[dl]_{w_0+w_1+w_2+w_3\choose w_0} 
\ar[d]|{{w_0+w_1+w_2+w_3\choose w_0+w_1}} 
\ar[dr]^{w_0+w_1+w_2+w_3\choose w_0+w_1+w_2} & \\
{\begin{tikzpicture}[scale=0.8]
\draw [thick] (1,0)--(2,0)--(3,0);
\node[arn_n] at (0,0) {};
\node[arn_n] at (1,0) {}; 
\node[arn_n] at (2,0) {}; 
\node[arn_n] at (3,0) {}; 
\end{tikzpicture}} 
\ar[d]|{-{w_1+w_2+w_3\choose w_1}} \ar[dr]|(.3){-{w_1+w_2+w_3\choose w_1+w_2}}  
& 
{\begin{tikzpicture}[scale=0.8]
\draw [thick] (0,0)--(1,0);
\draw [thick] (2,0)--(3,0);
\node[arn_n] at (0,0) {};
\node[arn_n] at (1,0) {}; 
\node[arn_n] at (2,0) {}; 
\node[arn_n] at (3,0) {}; 
\end{tikzpicture}} 
\ar[dl]^(.3){{w_0+w_1\choose w_0}} \ar[dr]_(.3){-{w_2+w_3\choose w_2}}  
&
{\begin{tikzpicture}[scale=0.8]
\draw [thick] (0,0)--(1,0)--(2,0);
\node[arn_n] at (0,0) {};
\node[arn_n] at (1,0) {}; 
\node[arn_n] at (2,0) {}; 
\node[arn_n] at (3,0) {}; 
\end{tikzpicture}} 
\ar[d]|{{w_0+w_1+w_2\choose w_0+w_1}} \ar[dl]|(.3){{w_0+w_1+w_2\choose w_0}}  
\\
{\begin{tikzpicture}[scale=0.8]
\draw [thick] (2,0)--(3,0);
\node[arn_n] at (0,0) {};
\node[arn_n] at (1,0) {}; 
\node[arn_n] at (2,0) {}; 
\node[arn_n] at (3,0) {}; 
\end{tikzpicture}} 
\ar[dr]_{w_2+w_3\choose w_2}
& 
{\begin{tikzpicture}[scale=0.8]
\draw [thick] (1,0)--(2,0);
\node[arn_n] at (0,0) {};
\node[arn_n] at (1,0) {}; 
\node[arn_n] at (2,0) {}; 
\node[arn_n] at (3,0) {}; 
\end{tikzpicture}} 
\ar[d]|{-{w_1+w_2\choose w_1}}
&
{\begin{tikzpicture}[scale=0.8]
\draw [thick] (0,0)--(1,0);
\node[arn_n] at (0,0) {};
\node[arn_n] at (1,0) {}; 
\node[arn_n] at (2,0) {}; 
\node[arn_n] at (3,0) {}; 
\end{tikzpicture}} 
\ar[dl]^{w_0+w_1\choose w_0}
\\
&
{\begin{tikzpicture}[scale=0.8]
\node[arn_n] at (0,0) {};
\node[arn_n] at (1,0) {}; 
\node[arn_n] at (2,0) {}; 
\node[arn_n] at (3,0) {}; 
\end{tikzpicture}}  & \\
}
\]
\caption{The complex $C_{\bullet}(w_0,w_1,w_2,w_3)$}
\label{fig:Cw-complex}
\end{figure}

\begin{problem}\label{prob:homology-Cw-complex}
  Describe the homology of the complex $C_{\bullet}(w_0,\cdots,w_d)$ over a field of characteristic $p$.
\end{problem}

To indicate the complexity of Problem~\ref{prob:homology-Cw-complex}, we discuss the answer in the case when $w_0=\cdots=w_d=1$, which is proved in \cite{rai-vdb}. To do so, we first introduce some notation. For a prime $p>0$, we enumerate non-negative integers $\equiv 0,1$ (mod $p$) as $0,1,p,p+1,2p,2p+1,\cdots$. If $m$ is in the list above then we write $|m|_p$ for its position, and call $|m|_p$ the \defi{$p$-index} of $m$:
\[ \text{if }m=pa+b,\text{ with }b\equiv 0,1\text{ (mod }p)\text{ and }0\leq b<p,\text{ then }|m|_p = 2a+b.\]
If $p=2$ then $|m|_p=m$, but when $p=3$ for instance, we get the following values of the $p$-index:
\[
\begin{array}{c|c|c|c|c|c|c|c|c|c|c|c|c|c}
m & 0 & 1 & 2 & 3 & 4 & 5 & 6 & 7 & 8 & 9 & 10 & 11 & \cdots \\ \hline
|m|_p & 0 & 1 &    &  2 & 3  &  & 4 & 5 &   & 6  &  7  &   & \cdots \\
\end{array}
\]
For a tuple $\a = (\a_0,\cdots,\a_k)$ of non-negative integers, with $\a_i\equiv 0,1\text{ mod }p$, we write
\[|\a|_p = \sum_{i=0}^k |\a_i|_p,\]
and we define
\[A_{p,d} = \left\{ \alpha=(\a_0,\cdots,\a_k) : \sum_{i=0}^k \a_i\cdot p^i = d,\ \a_i\geq 0,\ \a_i\equiv 0,1\text{ (mod $p$)}\right\}.\]
We have for instance $A_{p,d} = \emptyset$ when $d\not\equiv 0,1$ (mod $p$),
\begin{equation}\label{eq:some-Apd}
 A_{2,4} = \{(0,0,1),(0,2),(2,1),(4)\},\ A_{3,4} = \{(1,1),(4)\}, \text{ and}
\end{equation}
\begin{equation}\label{eq:some-p-indices}
 |(0,0,1)|_2 = 1,|(0,2)|_2 = 2,|(2,1)|_2 = 3,|(4)|_2 = 4,|(1,1)|_3 = 2,|(4)|_3 = 3.
\end{equation}

\begin{theorem}[\cite{rai-vdb}]\label{thm:hom-Call1} Let $C_{\bullet}=C_{\bullet}(w_0,\cdots,w_d)$, where $w_0=\cdots=w_d=1$, and suppose that $p=\op{char}(\kk)>0$. The following equality holds in $\bb{Z}[t]$:
  \[ \sum_{i\geq 0} \dim_{\kk} H_i(C_{\bullet}) \cdot t^i = \sum_{\a\in A_{p,d+1}} t^{d+1-|\a|_p}\]
\end{theorem}

We can then use \eqref{eq:some-Apd}, \eqref{eq:some-p-indices} and Theorem~\ref{thm:hom-Call1} to compute the homology of the complex $C_{\bullet}$ in \eqref{eq:C1111} when $p=2,3$. We have $d=3$ and therefore
\begin{equation}\label{eq:hom-C111} 
\sum_{i\geq 0} \dim_{\kk} H_i(C_{\bullet}) \cdot t^i = \begin{cases}
1+t+t^2+t^3 & \text{if }p=2,\\
t+t^2 &  \text{if }p=3.
\end{cases}
\end{equation}
Since $4\not\equiv 0,1$ (mod $p$) when $p\geq 5$, it follows that $C_{\bullet}$ is exact for such $p$.

It is difficult to imagine a general answer for arbitrary weights $w_0,\cdots,w_d$, but one could hope for a recursive description of the homology. This involves relating the homology of $C_{\bullet}(\ul{w})$ for different weight sequences $\ul{w}$. An easy consequence of Lucas' theorem shows for instance that
\begin{equation}\label{eq:Lucas} 
C_{\bullet}(w_0+p^r,w_1,\cdots,w_d) \simeq C_{\bullet}(w_0,w_1,\cdots,w_d) \quad\text{ if } p^r >w_1+\cdots+w_d.
\end{equation}
One can also check that for each $i$ there is a short exact sequence of complexes
\[ 0 \lra C_{\bullet}(w_0,\cdots,w_i) \oo C_{\bullet}(w_{i+1},\cdots,w_d) \lra C_{\bullet}(\ul{w}) \]
\[\lra C_{\bullet-1}(w_0,\cdots,w_i+w_{i+1},\cdots,w_d) \lra 0\]

We end this section with a more subtle, conjectural property of the complexes $C_{\bullet}(\ul{w})$. If we consider the ring of integer-valued polynomials 
\[ R = \bigoplus_{n\geq 0} \bb{Z}\cdot {x\choose n} \subset \bb{Q}[x],\quad\text{where}\quad {x\choose n} = \frac{x(x-1)\cdots(x-n+1)}{n!},\]
then we can define $C_{\bullet}(x,w_1,\cdots,w_n)$ as a complex of free $R$-modules for all non-negative integers $w_1,\cdots,w_n\geq 0$ (by choosing $w'(J,j_i)$ in \eqref{eq:def-pdeJ} to be the weight of the interval not containing $w_0=x$, so that the resulting binomial expression makes sense as a polynomial in $x$).
\begin{conjecture}\label{conj:involution-Cw}
 If we take $w_1=\cdots=w_d=1$ then we have an isomorphism of complexes of $R$-modules
 \[ C_{\bullet}(x,1,\cdots,1) \simeq C_{\bullet}(-x-2d,1,\cdots,1).\]
\end{conjecture}

Notice that we can specialize $x$ to an arbitrary integer and obtain a complex of $\bb{Z}$-modules, and then base change to any field $\kk$, thus obtaining complexes $C_{\bullet}(\ul{w})$ over an arbitrary field, where $w_0$ is allowed to be a negative number! Combining Conjecture~\ref{conj:involution-Cw} with \eqref{eq:Lucas}, one gets that in characteristic $p>0$
 \[ C_{\bullet}(w_0,1,\cdots,1) \simeq C_{\bullet}(p^r-w_0-2d,1,\cdots,1)\text{ for }p^r>w_0+2d,\]
 which we can verify independently. We do not know if there is a way to generalize the construction of $C_{\bullet}(\ul{w})$ to allow multiple weights to be negative, but this could likely lead to more symmetries among the complexes $C_{\bullet}(\ul{w})$.

\section{Trivial group action and stabilization of cohomology}\label{sec:stabilization}

If we identify $Fl_n = \GL_n/B$ where $B$ is the Borel subgroup of upper triangular matrices, then we can see that the $\GL_n$-equivariant line bundles on $Fl_n$ form a free abelian group of rank $n$. Algebraically, they are parametrized by one-dimensional representations of $B$ (or equivalently, characters of the maximal torus). Geometrically, they are obtained as follows (see also \cite[Section~9.2]{BCRV}): we write $\mc{Q}_i$ for the universal rank $i$ quotient sheaf on $Fl_n$, let $\mc{L}_i = \ker(\mc{Q}_i\onto\mc{Q}_{i-1})$ and define
\[ \mc{O}_{Fl_n}(\ll) = \mc{L}_1^{\ll_1}\oo\cdots\oo\mc{L}_n^{\ll_n} \quad\text{ for }\ll\in\bb{Z}^n.\]

The tautological surjection
\[ \kk^n \oo \mc{O}_{Fl_n} \onto \mc{Q}_{n-1}\]
induces an isomorphism on global sections, and we write $x_i$ for the global section corresponding to the $i$-th standard unit vector in $\kk^n$. We can then realize $Fl_{n-1}$ as the zero locus of the section $x_n\in H^0(Fl_n,\mc{Q}_{n-1})$, and the inclusion $Fl_{n-1}\subset Fl_n$ is $\GL_{n-1}$-equivariant. Moreover, for $\ll\in\bb{Z}^n$ we have
\begin{equation}\label{eq:restr-Olam}
 \mc{O}_{Fl_n}(\ll)_{|_{Fl_{n-1}}} = \mc{O}_{Fl_{n-1}}(\ol{\ll}) \text{ where }\ol{\ll} = (\ll_1,\cdots,\ll_{n-1}) \in \bb{Z}^{n-1}.
\end{equation}
Given a weight $\ll\in\bb{Z}^r$ we define
\[ \ll^{[n]} = (\ll_1,\cdots,\ll_r,0,\cdots,0) \in \bb{Z}^n \quad\text{ for }n\geq r,\]
and consider the maps
\[
\psi^j_n : H^j\left(Fl_n,\mc{O}_{Fl_n}(\ll^{[n]})\right) \lra H^j\left(Fl_{n-1},\mc{O}_{Fl_{n-1}}(\ll^{[n-1]})\right).
\]
induced on cohomology by the inclusion $Fl_{n-1}\subset Fl_n$ and the identification \eqref{eq:restr-Olam}. 

\begin{theorem}[\cite{rai-vdb}]\label{thm:stabilization}
 Suppose that $\ll\in\bb{Z}^r$ satisfies $|\ll|=\ll_1+\cdots+\ll_r\leq 0$.
 \begin{enumerate}
  \item The action of $\GL_n$ on $H^j\left(Fl_n,\mc{O}_{Fl_n}(\ll^{[n]})\right)$ is trivial for $n\gg 0$.
  \item The map $\psi^j_n$ is an isomorphism for $n\gg 0$.
  \item If $|\ll|<0$ then $H^j\left(Fl_n,\mc{O}_{Fl_n}(\ll^{[n]})\right)=0$ for $n\gg 0$.
 \end{enumerate}
\end{theorem}

When $\op{char}(\kk)=0$, the conclusions of Theorem~\ref{thm:stabilization} are easy applications of the Borel--Weil--Bott theorem, but in characteristic $p$ they are more subtle. Theorem~\ref{thm:stabilization} allows us to define for $|\ll|\leq 0$ the \defi{stable cohomology groups}
\[ H^j_{st}(\ll) = H^j\left(Fl_n,\mc{O}_{Fl_n}(\ll^{[n]})\right) \text{ for }n\gg 0\]
and formulate the stable versions of Problems~\ref{prob:BWB-char-p} and~\ref{prob:char-p-vanishing}.

\begin{problem}\label{prob:stab-coh-Olam}
 Compute the stable cohomology groups $H^j_{st}(\ll)$ over a field of characteristic $p>0$, and/or characterize their (non-)vanishing behavior.
\end{problem}

The significance of Theorem~\ref{thm:hom-Call1} can now be explained by its relation to stable cohomology. We show in \cite{rai-vdb} that
\begin{equation}\label{eq:stab-111}
 H^j_{st}(\ll) = H_{d+1-j}(C_{\bullet}(\ul{w}))  \text{ for }\ll = (-d-1,d+1)\in\bb{Z}^2,\ w_0=\cdots=w_d=1.
\end{equation}
In fact, we show more generally that if $w_0>0$, $w_1=\cdots=w_d$, and if we let
\begin{equation}\label{eq:hook-weight}
\ll = (-w_0-d,d+1,1,\cdots,1) \in \bb{Z}^{1+w_0}
\end{equation}
then we have
\[ H^j_{st}(\ll) = H_{d+w_0-j}(C_{\bullet}(\ul{w})).\]
It follows from \eqref{eq:hom-C111} that for $\ll=(-4,4)$ the only non-zero stable cohomology is
\[ H^1_{st}(\ll) = H^2_{st}(\ll) = H^3_{st}(\ll) = H^4_{st}(\ll) = \kk \quad\text{ when }\op{char}(\kk)=2,\text{ and }\]
\[ H^2_{st}(\ll) = H^3_{st}(\ll) = \kk \quad\text{ when }\op{char}(\kk)=3.\]
Perhaps one familiar example is obtained by taking $d=0$ in \eqref{eq:hook-weight}, where the cohomology is independent on the characteristic of $\kk$ and is given by
\[ H^{w_0}_{st}(\ll) = \kk \quad \text{ for }\ll=(-w_0,1,\cdots,1)\in\bb{Z}^{w_0+1}.\]
There are many natural questions that one can pose now, regarding the effective bounds for stabilization as well as efficient algorithmic methods for computing the stable cohomology.

\begin{question}\label{que:algo-for-stable-coh}
 Is it possible to construct arithmetic complexes such as the ones in Section~\ref{sec:arithmetic}, which compute the stable cohomology $H^j_{st}(\ll)$ for all $\ll$ and all characteristics?
\end{question}

We note here that the homology of the complexes $C_{\bullet}(\ul{w})$ for a general $\ul{w}$ does have an interpretation in terms of sheaf cohomology, but it is the cohomology of certain higher rank vector bundles on projective space, as explained in the next section.

\section{Polynomial functors of the cotangent sheaf on projective space}

If we let $\PP^{n-1} = \bb{P}(\kk^n)$ denote the $(n-1)$-dimensional projective space, then there is a natural forgetful map
\[ \pi:Fl_n \lra \PP^{n-1},\quad \pi(V_{\bullet}) = V_1,\]
and some of the cohomology groups from Section~\ref{sec:stabilization} can be realized as sheaf cohomology groups of vector bundles on $\PP^{n-1}$, as follows. Recall that a partition $\mu=(\mu_1,\mu_2,\cdots)$ is a non-decreasing sequence $\mu_1\geq\mu_2\geq\cdots\geq 0$. We write $\Omega$ for the \defi{cotangent sheaf} on $\PP^{n-1}$, and write $\bb{S}_{\mu}$ for the \defi{Schur functor} associated to $\mu$, noting that $\bb{S}_{\mu}\Omega=0$ if $\mu_n>0$ since $\op{rank}(\Omega)=n-1$. If $\mu_n=0$ then one can check using the projection formula (see \cite[Theorem~9.8.5, and Section~9.3]{BCRV}) that
\[ H^j\left(\PP^{n-1},\bb{S}_{\mu}\Omega(e)\right) = H^j\left(Fl_n,\mc{O}_{Fl_n}(\ll)\right) \text{ where }\ll = (e-|\mu|,\mu_1,\cdots,\mu_{n-1})\in\bb{Z}^n.\]
In the special case when $e=0$ and $\mu=(d,0,\cdots)$ we obtain
\[ H^j\left(\PP^{n-1},\Sym^d\Omega\right) = H^j\left(Fl_n,\mc{O}_{Fl_n}(\ll)\right) \text{ where }\ll=(-d,d,0,\cdots),\]
whose (stable) cohomology can be computed based on \eqref{eq:stab-111} and Theorem~\ref{thm:hom-Call1}.

In analogy with Theorem~\ref{thm:stabilization}, we can prove a general stabilization result, as well as a statement regarding the trivial $\GL_n$-action on cohomology.

\begin{theorem}[\cite{rai-vdb}]\label{thm:coh-omega}
 Suppose that $T$ is a polynomial functor of degree $d$.
 \begin{enumerate}
  \item The action of $\GL_n$ on $H^j\left(\PP^{n-1},T(\Omega)\right)$ is trivial for all $j<n-1$.
  \item If $n-1\geq d$ then the action of $\GL_n$ on $H^{n-1}\left(\PP^{n-1},T(\Omega)\right)$ is trivial.
  \item $H^j\left(\PP^{n-1},T(\Omega)\right)$ is independent on $n$ as long as $n-1\geq d$.
 \end{enumerate}
\end{theorem}

It is easy to see that the condition $n-1\geq d$ is necessary: for instance if $d=n=2$ and $T=\Sym^2$, then $\Omega=\omega$ is the canonical line bundle on $\PP^1$ and 
\begin{equation}\label{eq:H1-omega2}
H^1(\PP^1,\omega^{\oo 2}) = \op{D}^2(\kk^2) \oo \bw^2(\kk^2)^{\vee},
\end{equation}
where $\op{D}^d$ denotes the $d$-th divided power functor.

In light of Theorem~\ref{thm:coh-omega}, we can talk about \defi{stable cohomology} for $T(\Omega)$ and ask the analogues of Problem~\ref{prob:stab-coh-Olam} and Question~\ref{que:algo-for-stable-coh} in this context. In the case when $T$ is a Schur functor, we propose the following periodicity conjecture, which we can verify when $l=1$ or $\mu_2\leq 2$.
\begin{conjecture}\label{conj:periodicity}
 Consider a partition $\mu=(\mu_1,\cdots,\mu_l)$ with $\mu_l>0$, and write 
 \[\mu[q] = (\mu_1,\cdots,\mu_l,1,\cdots,1)\] 
 for the partition obtained from $\mu$ by appending $q$ parts of size $1$. If $q=p^r>|\mu|-l$ and $\op{char}(\kk)=p$ then 
 \[H^j(\PP^{n-1},\bb{S}_{\mu}\Omega) \simeq H^{j+q}(\PP^{n-1},\bb{S}_{\mu[q]}\Omega).\]
\end{conjecture}

Schur functors can be associated not only to partitions, but also to \defi{skew shapes}~$\ll/\mu$ \cite[Section~II.1]{ABW}, where $\ll,\mu$ are partitions with $\mu_i\leq\ll_i$ for all $i$. We picture a shape $\ll/\mu$ by subtracting the Young diagram of $\mu$ from that of $\ll$, and say that $\ll/\mu$ is a \defi{ribbon} if it contains no $2\times 2$ square. For instance if $\ll=(7,4,3,3,3)$ and $\mu=(3,2,2,2)$ then we get the ribbon skew-shape
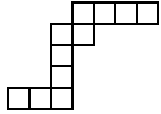
\[
\ytableausetup{boxsize=0.7em}
\begin{ytableau} 
\none& \none& \none & & & &\\ 
 \none& \none&  & \\  
 \none& \none&   \\  
 \none& \none&    \\  
   & &    \\  
\end{ytableau}
\]
which we can also record using the sequence of column sizes $\ul{w} = (1,1,4,2,1,1,1)$. Note that if $w_0\geq 1$ and $w_1=\cdots=w_d=1$ then the corresponding skew-shape is a hook shape, and hence an actual partition. The significance of the complexes $C_{\bullet}(\ul{w})$ from Section~\ref{sec:arithmetic} is that they compute the stable cohomology of $\bb{S}_{\ll/\mu}\Omega$, where $\ll/\mu$ is the skew-shape corresponding to $\ul{w}$. The statement of Conjecture~\ref{conj:periodicity} is then natural in light of \eqref{eq:Lucas}.

The relationship between the cohomology of $\bb{S}_{\ll/\mu}\Omega$ and the homology of $C_{\bullet}(\ul{w})$ is based on the hypercohomology spectral sequence associated with the resolution of $\bb{S}_{\ll/\mu}$ by tensor product of exterior power functors. The question of resolving Schur functors by tensor products of exterior (or divided) powers has a long history, that played a significant part in David Buchsbaum's work \cite{AB1,AB2,BR1,BR2,BR3,Buc}. Akin and Buchsbaum proved already in the 80s that such resolutions exist, but the quest to find an optimal one continues on: perhaps the best general results are due to Santana and Yudin \cite{san-yud}, but their constructions give resolutions that are much longer than the global dimension of the relevant category of polynomial functors \cite[Theorem~2]{totaro}.

\begin{problem}\label{prob:res-by-ext-powers}
 Find the (an) \emph{optimal} resolution of a Schur functor $\bb{S}_{\ll}$ (or $\bb{S}_{\ll/\mu}$) by tensor products of exterior power functors.
\end{problem}

The discussion above simplifies quite a bit if we consider Weyl functors $\bb{W}_{\mu}$ instead of Schur functors -- a quick way to define them is by letting 
\[\bb{W}_{\mu}(V) = (\bb{S}_{\mu}(V^{\vee}))^{\vee}\]
for a finite dimensional $\kk$-vector space $V$. If $\mu_1\geq 2$ and $n>|\mu|$ then one can show
\[ H^j\left(\PP^{n-1},\bb{W}_{\mu}\Omega\right) = 0\text{ for all }j.\]
If $\mu_1=1$ and $d=|\mu|$ then we have $\bb{W}_{\mu}\Omega = \bw^d \Omega = \Omega^d$ is the sheaf of differential $d$-forms, in which case the only non-vanishing cohomology is $H^d\left(\PP^{n-1},\Omega^d\right) = \kk$. To make the question of computing cohomology more interesting, one should then consider twists by a line bundle: we formulate the following general problem, whose difficulty is likely comparable to that of Problems~\ref{prob:BWB-char-p} and~\ref{prob:char-p-vanishing}.

\begin{problem}\label{prob:coh-table-TOmega}
 Given a polynomial functor $T$, describe (and characterize the vanishing behavior of) the cohomology of $T(\Omega)\oo\mc{O}_{\PP^{n-1}}(e)$ for all $e\in\bb{Z}$.
\end{problem}

\section{The incidence correspondence}

Problem~\ref{prob:coh-table-TOmega} is already open in the case when $T=\op{D}^d$ is a \defi{divided power functor} (the Weyl functor $\bb{W}_{\mu}$ for $\mu=(d,0,\cdots)$). The vanishing behavior of the cohomology of $\op{D}^d\Omega(e)$ is characterized in \cite{Gao-Raicu}, but we do not have a general description of the cohomology groups. To provide some geometric context, we consider the variety $X$ parametrizing partial flags
\[ 0 \subset V_1 \subset V_{n-1} \subset \kk^n,\]
which in the case $n=3$ agrees with $Fl_3$. If we write $\PP=\PP^{n-1}$, and let $\PP^{\vee}$ denote the dual projective space, we can reinterpret $X$ as the \defi{incidence correspondence}
\[ X = \{(p,H) : p\in H\} \subset \PP \times \PP^{\vee},\]
which is a hypersurface cut out by the bilinear form
\begin{equation}\label{eq:def-omega}
 \omega = x_1y_1+\cdots+x_ny_n.
\end{equation}
The line bundles on $X$ arise by restriction from $\PP \times \PP^{\vee}$, and we write
\[ \mc{O}_X(a,b) = \mc{O}_{\PP \times \PP^{\vee}}(a,b)_{|_X}\quad\text{ for }(a,b)\in\bb{Z}^2.\]
Along the first projection $\varphi:X\lra\PP$, we can identify $X$ with the projective bundle $\bb{P}_{\PP}(\Omega) = \ul{\op{Proj}}\left(\Sym(\Omega^{\vee})\right)$ which after some manipulations yields
\[ H^j\left(\PP, \op{D}^d\Omega \oo\mc{O}_{\PP}(e)\right) \simeq H^{j+n-2}\left(X,\mc{O}_X(e-d+1,-d-n+1)\right) \oo \left(\bw^n \kk^n\right)^{\vee},\]
where the last tensor factor is there to make the isomorphism $\GL_n$-equivariant. In terms of the forgetful map 
\[ \psi:Fl_n \lra X,\quad \psi(V_{\bullet}) = \left(V_1 \subset V_{n-1}\right),\]
we have that $\psi^*\left(\mc{O}_X(a,b)\right) = \mc{O}_{Fl_n}(a,0,\cdots,0,-b)$ and by the projection formula
\[ H^j\left(X,\mc{O}_X(a,b)\right) = H^j\left(Fl_n,\mc{O}_{Fl_n}(a,0,\cdots,0,-b)\right),\]
showing that Problem~\ref{prob:coh-table-TOmega} for $T=\op{D}^d$ is just a special case of Problems~\ref{prob:BWB-char-p},~\ref{prob:char-p-vanishing}.

We find it more convenient to rephrase the cohomology calculations in terms of the \defi{tautological sheaf} $\mc{R}=\Omega(1)$ (as in \cite{Gao-Raicu}), and note that we have
\[ 
\begin{aligned}
H^j\left(\PP, \op{D}^d\mc{R} \oo\mc{O}_{\PP}(e)\right) &= H^j\left(\PP, \op{D}^d\Omega \oo\mc{O}_{\PP}(e+d)\right) \\
&\simeq H^{j+n-2}\left(X,\mc{O}_X(e+1,-d-n+1)\right)\oo \left(\bw^n \kk^n\right)^{\vee}.
\end{aligned}
\]
For yet another equivalent interpretation of these cohomology groups, we let 
\begin{equation}\label{eq:defR}
 R = \kk[x_1,\cdots,x_n,y_1,\cdots,y_n],
\end{equation}
and consider the (infinitely generated) local cohomology module
\[ M = H^n_{(y_1,\cdots,y_n)}(R).\]
We define a bi-grading on $M$ by letting
\[ M_{d,e} = \bigoplus_{\substack{a_i,b_j\geq 0 \\ |\ul{a}|=d,\ |\ul{b}|=e}} \kk\cdot\frac{x_1^{b_1}\cdots x_n^{b_n}}{y_1^{1+a_1}\cdots y_n^{1+a_n}}.\]
For $\omega$ as in \eqref{eq:def-omega}, $\mc{L}=\mc{O}_X(e+1,-d-n+1)$, $e\geq -1$, we have an exact sequence
\[ 0 \lra H^{n-2}\left(X,\mc{L}\right) \lra M_{d,e} \overset{\cdot\omega}{\lra} M_{d-1,e+1} \lra H^{n-1}\left(X,\mc{L}\right)\lra 0\]
and $H^j\left(X,\mc{L}\right)=0$ for $j\not\in\{n-2,n-1\}$ (when $e\leq -2$ the cohomology can only be non-zero for $j=2n-3$, and its dimension is the same in all characteristics). For a concrete example that illustrates the dependence on characteristic, take $n=3$, $d=2$, $e=1$ and consider the element
\begin{equation}\label{eq:coh-class-f}
f = \frac{x_1}{y_1y_2^2y_3^2}+\frac{x_2}{y_1^2y_2y_3^2}+\frac{x_3}{y_1^2y_2^2y_3} = \frac{\omega}{y_1^2y_2^2y_3^2} \in M_{2,1}.
\end{equation}
We have that
\[ f\cdot\omega = \frac{\omega^2}{y_1^2y_2^2y_3^2} = 2\cdot\left(\frac{x_1x_2y_1y_2+x_1x_3y_1y_3+x_2x_3y_2y_3}{y_1^2y_2^2y_3^2}\right)\]
which is zero if and only if $\op{char}(\kk)=2$. In that case we have that $H^1\left(X,\mc{L}\right) = \kk\cdot f$ is one dimensional.

As opposed to the examples in previous sections, the cohomology groups we consider here usually have a non-trivial $\GL_n$-action, which should play a role in describing their structure. A first approximation of this comes from the action of the maximal torus $(\kk^{\times})^n$ of diagonal matrices. This is encoded by the notion of a \defi{character}, or equivalently by a $\bb{Z}^n$-grading: in the preceding paragraphs, the natural $\bb{Z}^n$-grading comes from assigning $\op{deg}(x_i)=\vec{e}_i$ (the $i$-th standard unit vector), and $\op{deg}(y_i)=-\vec{e}_i$. We review the notion of characters next, followed by a discussion of some conjectural character formulas.

\subsection{Characters}

Representations of the algebraic torus $(\kk^{\times})^n$ are equivalent to $\bb{Z}^n$-graded vector spaces \cite[Lemma~9.7.9]{BCRV}, \cite[Section~2.11]{Jantzen}, and for any such (finite dimensional) vector space $W$, its \defi{character} is
\[[W] := \sum_{(i_1,\cdots,i_n)\in\bb{Z}^n} \dim(W_{(i_1,\cdots,i_n)})\cdot t_1^{i_1}\cdots t_n^{i_n} \in \bb{Z}[t_1^{\pm 1},\cdots,t_n^{\pm 1}].\]
If $W$ is a $\GL_n$-representation, then $[W]$ is invariant under the action of the symmetric group $\mf{S}_n$ by permutations of $t_1,\cdots,t_n$. Examples include the \defi{complete symmetric polynomials $h_d$} and \defi{Schur polynomials $s_{(a,b)}$}, defined by
\[ h_d = [\Sym^d(\kk^n)] = \sum_{\substack{i_1+\cdots+i_n = d \\ i_j\geq 0}}t_1^{i_1}\cdots t_n^{i_n},\quad s_{(a,b)} = [\bb{S}_{(a,b)}(\kk^n)] = h_a\cdot h_b-h_{a+1}\cdot h_{b-1}.\]
Notice that although $\Sym^d(\kk^n)$ and $\op{D}^d(\kk^n)$ are not isomorphic $\GL_n$-representations in general, they are isomorphic as $(\kk^{\times})^n$-representations and therefore they have the same character. We have
\[ M_{d,e}\simeq \op{D}^d\left(\kk^n\right)\oo\Sym^e\left(\kk^n\right)\oo\left(\bw^n\kk^n\right),\quad[M_{d,e}] = h_d\cdot h_e\cdot(t_1\cdots t_n).\]
Since the element $f$ in \eqref{eq:coh-class-f} has multidegree $(2,2,2)$, we get for $n=3$, $\op{char}(\kk)=2$
\[\left[ H^1\left(X, \mc{O}_X(2,-4)\right) \right] =t_1^2t_2^2t_3^2 \quad\text{and}\quad \left[ H^0\left(\PP, \op{D}^2\mc{R}(1)\right) \right] = t_1t_2t_3.\]
For a non-polynomial example, if $W$ is the representation in \eqref{eq:H1-omega2} then
\[ [W] = t_1^2 t_2^{-1} + 1 + t_1^{-1}t_2^2.\]

To describe further cohomology characters, it is useful to consider \defi{$q$-truncated complete symmetric polynomials $h_d^{(q)}$} and \defi{$q$-truncated Schur polynomials $s^{(q)}_{(a,b)}$}, defined by (see also \cite{DW}, \cite{Walker})
\[
h^{(q)}_d = \sum_{\substack{i_1+\cdots+i_n = d \\ 0\leq i_j<q}}t_1^{i_1}\cdots t_n^{i_n},\qquad s^{(q)}_{(a,b)} = h^{(q)}_a\cdot h^{(q)}_b-h^{(q)}_{a+1}\cdot h^{(q)}_{b-1}.
\]
as well as the \defi{$q$-Frobenius endomorphisms} of the character ring
\[F^q:\bb{Z}[t_1^{\pm 1},\cdots,t_n^{\pm 1}]\lra \bb{Z}[t_1^{\pm 1},\cdots,t_n^{\pm 1}],\quad\text{ defined by }t_i\mapsto t_i^q.\]

\subsection{Weights close to the origin} It is a general phenomenon that the cohomology characters of $\mc{O}_{Fl_n}(\ll)$ behave well (as in characteristic zero) when $\ll$ is close to the origin, and get increasingly more complicated the further away one gets from the origin. For the divided powers of $\mc{R}$ (or $\Omega$) this translates into the fact that the cohomology of $\op{D}^d\mc{R}(e)$ in characteristic $p$ is the same as that in characteristic zero when $0\leq d<p$. Passing to the next order of magnitude ($p\leq d<p^2$) we have the following conjectural description of the cohomology characters (see \cite{Gao-Raicu} for the explanation of the reduction to the case $e\geq d-1$).

\begin{conjecture}\label{conj:small-weights}
 Suppose that $tp\leq d<(t+1)p$, where $1\leq t<p$. We have for all $e\geq d-1$
 \[ \left[H^1\left(\PP^{n-1},D^d\mc{R}(e)\right)\right] = \sum_{\substack{1\leq b\leq a\leq t \\ 0\leq j\leq a-b}} F^p(s_{(a-b,j)})\cdot s^{(p)}_{(e+(b-j)p,d-ap)}.\]
\end{conjecture}

When $t=1$, Conjecture~\ref{conj:small-weights} simplifies to the assertion that 
\begin{equation}\label{eq:H1-d<2p}
\left[H^1\left(\PP^{n-1},\op{D}^d\mc{R}(e)\right)\right] = s^{(p)}_{(e+p,d-p)},
\end{equation}
which is proved in \cite[Theorem~1.6]{Gao-Raicu}. When $n=3$, Conjecture~\ref{conj:small-weights} can be proved using \cite[Theorem~1.7]{Gao-Raicu}. In general, we expect that a character formula for the cohomology of $\op{D}^d\mc{R}(e)$ will rely in an essential way on truncated Schur polynomials, and we formulate a precise statement when $\op{char}(\kk)=2$ below.

\subsection{Nim polynomials} In addition to truncated Schur polynomials, the other ingredient in the description of cohomology characters when $\op{char}(\kk)=2$ appears to come from the symmetric polynomials defining winning positions in the game of Nim (see \cite{Bouton} for the basic theory), which we introduce next.

We write $d=(d_k\cdots d_0)_2$ for the $2$-adic expansion of a non-negative integer $d$:
\[ d= \sum_{i=0}^k d_i\cdot 2^i,\text{ with }d_i\in\{0,1\}\text{ for all }i.\]
The \defi{Nim-sum} (or \defi{bitwise xor}) $a\oplus b$ is defined by performing addition modulo $2$ to each of the digits in the $2$-adic expansion of non-negative integers $a,b$:
\[ a\oplus b = c\text{ if and only if }a_i + b_i \equiv c_i \text{ mod }2\text{ for all }i,\]
where $a_i,b_i,c_i$ denote the digits of $a,b,c$ in the $2$-adic expansion. We define the $n$-variate \defi{Nim symmetric polynomials} via
\[
 \mc{N}_m = \sum_{\substack{i_1+\cdots+i_n=2m \\ i_1\oplus i_2\oplus\cdots\oplus i_n = 0}} t_1^{i_1} t_2^{i_2}\cdots t_n^{i_n}.
\]
We have for instance that $\mc{N}_0=1$ and
\[ \mc{N}_1 = \sum_{1\leq i<j\leq n} t_i t_j = s_{(1,1)}.\]

\begin{conjecture}\label{conj:char2-characters}
  If $\op{char}(\kk)=2$ and $e\geq d-1$ then
 \[ \left[H^1\left(\PP^{n-1},D^d\mc{R}(e)\right)\right] = \sum_{(q,m)\in\Lambda_d} F^{2q}(\mc{N}_m)\cdot s^{(q)}_{(e-(2m-1)q,d-(2m+1)q)},\]
 where $\Lambda_d = \left\{(q,m) | q=2^r\text{ for some }r\geq 1,\ m\geq 0,\text{ and }(2m+1)q\leq d\right\}$.
\end{conjecture}

When $d=2,3$ we have $\Lambda_d = \{(2,0)\}$ and Conjecture~\ref{conj:char2-characters} predicts
\[ \left[H^1\left(\PP^{n-1},D^d\mc{R}(e)\right)\right] = s^{(2)}_{(e+2,d-2)}\]
which is a special case of \eqref{eq:H1-d<2p}. When $d=6$ we have $\Lambda_d =  \{(4,0), (2,0), (2,1)\}$ and the conjecture asserts
\[ \left[H^1\left(\PP^{n-1},D^6\mc{R}(e)\right)\right] = s^{(4)}_{(e+4,2)} + s^{(2)}_{(e+2,4)} + F^4(\mc{N}_1) \cdot s^{(2)}_{(e-2,0)} \quad\text{ for }e\geq 5.\]
We have checked computationally that the dimensions of the corresponding representations (obtained by setting $t_1=\cdots=t_n=1$) agree for $n\leq 10$.

Conjecture~\ref{conj:char2-characters} appears in \cite[Chapter~4]{Gao-thesis} where it is proved for $n=3,4$, and it is also discussed in an equivalent form for $n=3$ in \cite[Theorem~1.9]{Gao-Raicu}. Computational data suggest character formulas with a similar flavor in other characteristics, but we do not know how to construct an appropriate replacement for the Nim polynomials $\mc{N}_m$ in general. 

\section{Relation to determinantal ideals}\label{sec:determinantal}

Schur functors are intimately related to the study of determinantal rings and ideals, particularly through the standard monomial theory for the ring $R$ of polynomial functions on the space of $m\times n$ matrices \cite{DCEP}, \cite[Chapter~3]{BCRV}. In characteristic zero it is now well-understood how to classify the ideals $I\subset R$ which are $\GL_m\times\GL_n$-invariant (for the action by row and column operations on the matrix entries) as well as how to compute basic homological invariants such as their Hilbert function. In characteristic $p>0$, even the case $m=1$ of the classification problem ($\GL_n$-invariant ideals in $\kk[x_1,\cdots,x_n]$) is non-trivial, and requires understanding the submodule structure for the $\GL_n$-representations $\Sym^d(\kk^n)$, which was worked out by Doty \cite{Doty}. 

In what follows we focus on the case of $2\times n$ matrices. We let $R$ as in \eqref{eq:defR}, and let $I$ denote the ideal of $2\times 2$ minors of the generic matrix
\[ \begin{bmatrix} x_1 & x_2 & \cdots & x_n \\ y_1 & y_2 & \cdots & y_n \end{bmatrix}.\]
There is a bi-grading on $R$ given by $\deg(x_i)=(1,0)$, $\deg(y_j)=(0,1)$ which makes the ideal $I$ bihomogeneous. The $\GL_n$-action (by column operations) makes~$R$ and~$I$ into $\GL_n$-representations, so we can analyze the characters of their bigraded components. We have for $a\geq b$
\[ \left[R_{(a,b)} \right] = \left[\Sym^a(\kk^n) \oo \Sym^b(\kk^n) \right] = s_a \cdot s_b = \sum_{i=0}^b s_{(a+b-i,i)},\]
where the last equality comes from Pieri's rule. This identity at the level of characters can be explained using the $I$-adic filtration on $R$: we have
\begin{equation}\label{eq:Iadic-deg-ab}
 \left[(I^i/I^{i+1})_{(a,b)}\right] = s_{(a+b-i,i)}\text{ for } i=0,\cdots,b,
\end{equation}
and $(I^i/I^{i+1})_{(a,b)} = 0$ for $i>b$. We can refine this further using tableaux combinatorics, focusing on the case $i=b$ for simplicity. We recall that a Young tableau 
 \begin{equation}\label{eq:general-2row-T}
 \ytableausetup{centertableaux,boxsize=1.5em}
T = \begin{ytableau}
u_1 & u_2 & \cdots & u_b & \cdots & u_a \\
v_1 & v_2 & \cdots & v_b
\end{ytableau}
\end{equation}
of shape $(a,b)$ is \defi{semi-standard} if 
\[u_1\leq u_2 \leq \cdots \leq u_a,\quad v_1\leq v_2\leq\cdots \leq v_b,\quad u_i<v_i\text{ for }1\leq i\leq b.\]
If we write $Tab_{n}(a,b)$ for the set of semi-standard tableaux of shape $(a,b)$ with entries in $\{1,\cdots,n\}$, and let
\[ t^T = t_{u_1}\cdots t_{u_a}\cdot t_{v_1} \cdots t_{v_b} \quad\text{ for }T\text{ as in }\eqref{eq:general-2row-T},\]
then the Schur polynomial $s_{(a,b)}$ can be computed as
\begin{equation}\label{eq:schur-from-tabs}
 s_{(a,b)} = \sum_{T\in Tab_{n}(a,b)} t^T.
\end{equation}
This suggests the existence of a basis of $(I^b)_{(a,b)} = (I^b/I^{b+1})_{(a,b)}$ indexed by $Tab_{n}(a,b)$, and indeed one such basis is constructed by letting
\[ G_T = \left(\prod_{i=1}^b (x_{u_i} y_{v_i} - x_{v_i} y_{u_i})\right) \cdot x_{u_{b+1}} \cdots x_{u_a}.\] 
The linear independence of the polynomials $G_T$ follows for instance from the fact that the leading term of $G_T$ with respect to the graded lexicographic order on $R$ (where $x_1>\cdots>x_n>y_1>\cdots>y_n$) is given by
\begin{equation}\label{eq:def-MT}
 M_T = x_{u_1}\cdots x_{u_a}\cdot y_{v_1}\cdots y_{v_b}.
\end{equation}

If one wants the analogues of this classical theory modulo Frobenius, then it is natural to consider (when $\op{char}(\kk)=p>0$) the quotient ring
\[ \ol{R} = R / \langle x_1^p,\cdots,x_n^p,y_1^p,\cdots,y_n^p\rangle,\]
and the corresponding ideal $\ol{I}$. The following analogue of \eqref{eq:Iadic-deg-ab} is already non-trivial.

\begin{conjecture}\label{conj:Iadic-filtration}
 If $a-b\geq p-1$ then 
 \[ \left[\left(\ol{I}^i/\ol{I}^{i+1}\right)_{(a,b)}\right] = s^{(p)}_{(a+b-i,i)}\text{ for } i=0,\cdots,b.\]
\end{conjecture}

We don't have a good explanation for the hypothesis $a-b\geq p-1$ (see also Conjecture~\ref{conj:compos-TpaTpb}), but it is easy to see by example that it is necessary: for instance when $a=b=1$, $i=0$, and $p=2$, we have $\left[(\ol{R}/\ol{I})_{(1,1)}\right] = h_2$ but $s^{(2)}_{(2,0)} = h^{(2)}_2$ is the second elementary symmetric polynomial. 

As before, we can look for standard bases that match up with the (conjectural) character formulas. The following appears to be the appropriate replacement of semi-standard tableaux.

\begin{definition}\label{def:weakly-pstd-tab}
We say that a tableau $T$ is \defi{$p$-semi-standard} if
\begin{enumerate}
    \item $T$ is weakly increasing along rows and columns: 
    \[u_1\leq\cdots\leq u_a,\quad v_1\leq\cdots \leq v_b,\quad\text{and}\quad u_i\leq v_i\text{ for }1\leq i\leq b.\]
    \item Any constant sequence within each row has length at most $(p-1)$: 
    \[u_i<u_{i+p}\quad\text{ and }\quad v_j<v_{j+p}\quad \text{for all }i,j.\] 
    \item If $u_j=v_j$, consider the unique indices $r\geq j$ and $s\leq j$ so that
    \[ u_j = \cdots = u_r < u_{r+1},\quad v_{s-1}<v_s=\cdots=v_j.\]
    We have $(r-j+1) + (j-s+1) \geq p$, that is, the total number of entries equal to $u_j=v_j$ to the right of (and including) $u_j$ and to the left of (and including) $v_j$ is at least $p$.
\end{enumerate}
We write $Tab^{(p)}_n(a,b)$ for the set of $p$-semi-standard tableaux of shape $(a,b)$ with entries in $\{1,\cdots,n\}$.
\end{definition}

Note that $T$ is $2$-semi-standard if and only if the transposed tableau to $T$ is semi-standard in the usual sense. When $p=3$, it is not hard to check that
\begin{equation}\label{eq:Tab3-21}
\ytableausetup{centertableaux,boxsize=1em}
Tab^{(3)}_n(2,1) = Tab_n(2,1) \bigcup \left\{\begin{ytableau}
i & i \\
i \\
\end{ytableau} \left| i=1,\cdots,n \right.\right\}
\end{equation}
One can show (without restrictions on $a,b$) the following analogue of \eqref{eq:schur-from-tabs}:
 \[ s^{(p)}_{(a,b)} = \sum_{T\in Tab^{(p)}_n(a,b)} t^T,\]
so \eqref{eq:Tab3-21} translates into the identity
\[ s^{(3)}_{(2,1)} = s_{(2,1)} + \sum_{i=1}^n t_i^3.\]
Note that truncated Schur polynomials may be \emph{larger} then their classical analogues, which is perhaps counterintuitive. It is also the case that for a general partition, a truncated Schur polynomial does not have to be the character of a $\GL_n$-representation (but will be a virtual character). Conjecture~\ref{conj:Iadic-filtration} asserts then that $s^{(p)}_{(a,b)}$ is indeed a character whenever $a-b\geq p-1$, and predicts a specific $\GL_n$-representation realizing it. A natural approach to Conjecture~\ref{conj:Iadic-filtration} in the case $i=b$ (which is in fact the main case) comes from the following.

\begin{conjecture}\label{conj:lead-terms}
 If $a-b\geq p-1$ then for every tableau $T\in Tab^{(p)}_n(a,b)$, there exists an element $H_T\in \ol{I}^b$ whose leading monomial is $M_T$.
\end{conjecture}

If we write $F^p\subset\Sym^p$ for the Frobenius power functor, and consider the \defi{truncated symmetric power} functors
\[T_p\Sym^d(\kk^n) = \op{coker}\left( F^p(\kk^n) \oo \Sym^{d-p}(\kk^n)  \lra \Sym^{d}(\kk^n)\right)\]
then one can view the discussion above as an attempt to understand 
\[\ol{R}_{(a,b)} = T_p\Sym^a(\kk^n) \oo T_p\Sym^b(\kk^n),\]
the tensor product of two truncated powers. This is a representation of $\GL_n$ (with character $h^{(p)}_a\cdot h^{(p)}_b$) and as such it has a filtration with composition factors given by simple modules $L(\ll)$ of highest weight $\ll$, for certain dominant weights $\ll$. 

\begin{conjecture}\label{conj:compos-TpaTpb}
 If $a-b\geq p-1$ then every composition factor $L(\ll)$ of $\ol{R}_{(a,b)}$ is \defi{$p$-restricted}, that is, $\ll_i - \ll_{i+1} \leq p-1$  for all $i$.
\end{conjecture}

The hypothesis $a-b\geq p-1$ is again necessary: if $a=p-1$, $b=1$ then $\ol{R}_{(a,b)}=\Sym^{p-1}(\kk^n) \oo \kk^n$ surjects onto $\Sym^p(\kk^n)$ which contains $F^p(\kk^n)$ as a submodule. Since $F^p(\kk^n)=L(\ll)$ for $\ll=(p,0,\cdots,0)$, this is not $p$-restricted. The question of understanding composition factors for tensor products of (any number of) truncated symmetric powers was considered in \cite{Don-Ger}, but the results there don't seem to be enough to settle Conjecture~\ref{conj:compos-TpaTpb}. There is of course no reason to restrict oneself to matrices with two rows, and we end with the hope that one day there will be a \defi{standard monomial theory modulo Frobenius}.

\section*{Acknowledgements}
We thank Alessio Sammartano for comments, corrections, and a careful reading of the manuscript.

\bibliographystyle{amsplain}

\end{document}